\documentclass[leqno,12pt]{amsart}
\usepackage{latexsym}
\usepackage{mathrsfs}
\usepackage{amsmath}
\usepackage{amssymb}
\usepackage[bookmarks]{hyperref}
\usepackage{pstricks,pst-plot}

\newtheorem{theorem}{Theorem}[section]
\newtheorem{lemma}[theorem]{Lemma}
\newtheorem{corollary}[theorem]{Corollary}
\newtheorem{proposition}[theorem]{Proposition}

\newtheorem{conjecture}[theorem]{Conjecture}
\newtheorem{observation}[theorem]{Observation}

\theoremstyle{definition}
\newtheorem{remark}[theorem]{Remark}
\newtheorem{definition}[theorem]{Definition}

\newtheorem{question}[theorem]{Question}

\newcommand{\C}{\mathbb{C}}

\newcommand{\Z}{\mathbb{Z}}

\newcommand{\R}{\mathbb{R}}

\newcommand{\bthm}{\begin{theorem}}
\newcommand{\ethm}{\end{theorem}}
\newcommand{\blem}{\begin{lemma}}
\newcommand{\elem}{\end{lemma}}
\newcommand{\bcor}{\begin{corollary}}
\newcommand{\ecor}{\end{corollary}}
\newcommand{\bprop}{\begin{proposition}}
\newcommand{\eprop}{\end{proposition}}
\newcommand{\bdefn}{\begin{definition}}
\newcommand{\edefn}{\end{definition}}
\newcommand{\bpf}{\begin{proof}}
\newcommand{\epf}{\end{proof}}

\def\vep {\varepsilon}
\def \sm {\setminus}

\def\itemskip {\vskip 3pt plus 2 pt minus 1 pt}
\def\sf{{\mathscr F}}

\def\xcf{X\times\C^\sf}

\def\aa{A_\alpha}

\newcommand{\ra}{\rightarrow}

\newcommand{\ol}{\overline}

\def\Dk{\Delta_k}
\def\unionDk{\bigcup_{k=1}^\infty \Dk}

\def\sd{{\mathscr D}}

\def\Rsphere{\C_\infty}

\begin{document}
\title{A sharper Swiss cheese}

\author{Alexander J. Izzo}
\address{Department of Mathematics and Statistics, Bowling Green State University, Bowling Green, OH 43403}
\email{aizzo@bgsu.edu}
\thanks{The author was partially supported by NSF Grant DMS-1856010.}


\begin{abstract}
It is shown$\vphantom{\widehat{\widehat{\widehat{\widehat{\widehat{\widehat{\widehat X}}}}}}}$
that there exists a compact planar set $K$ such that the uniform algebra $R(K)$ is nontrivial and strongly regular.  This settles an issue raised by Donald Wilken 55 years ago.  It is shown that the set $K$ can be chosen such that, in addition, $R(K)$ is not weakly amenable.  It is also shown that there exists a uniform algebra that has bounded relative units but is not weakly amenable.  These results answer questions raised by Joel Feinstein and Matthew Heath 17 years ago.  A key ingredient in our proofs is a bound we establish on the functions introduced by Thomas K\"orner to simplify Robert McKissick's construction of a nontrivial normal uniform algebra.
\end{abstract}

\maketitle

\vskip -2.56 true in
\centerline{\footnotesize\it Dedicated to the memory of Harold Garth Dales} 
\vskip 2.56 truein

%
%
%
%

\section{Introduction}

In this paper we answer several questions in the literature regarding strong regularity and weak amenability of uniform algebras.  (These terms are defined later in this introduction.  For definitions of other terminology and notation used in this introduction see Section~\ref{notation}.)
Our main goal is to prove the following.

\bthm\label{main-theorem}
There exists a compact set $K$ in the complex plane such that $R(K)$ is a nontrivial strongly regular uniform algebra.
\ethm

Here, as usual, $R(K)$ denotes the uniform closure on $K$ of the holomorphic rational functions with poles off $K$.

A uniform algebra $A$ on a compact space $X$ is \emph{strongly regular} at a point $x\in X$ if $\ol{J_x}=M_x$, where the ideals $M_x$ and $J_x$ are defined by the equations
\begin{align} 
M_x &= \{\, f\in A:  f(x)=0\,\}  \nonumber  \\ 
\intertext{and} 
J_x&= \{ \, f\in A: f^{-1}(0)\ \hbox{contains a neighborhood of $x$ in $X$}\}. \nonumber 
\end{align}
The uniform algebra $A$ is \emph{strongly regular} if $A$ is strongly regular at every point of $X$.
It was shown by Wilken that every strongly regular uniform algebra is normal \cite[Corollary~1]{Wilken2}.

The issue of whether $R(K)$, for $K$ a compact planar set, can be nontrivial and strongly regular (answered by the above theorem) was raised by 
Donald Wilken 55 years ago \cite{Wilken2}.  The question was reiterated 32 years ago by Joel Feinstein in the paper in which he constructed the first example of a nontrivial strongly regular uniform algebra~\cite{F1}, and the question was reiterated again by Feinstein and Matthew Heath in \cite{FH}.

In fact, our main result gives more detailed information than the above theorem.  Here, and throughout the paper, we denote the open unit disc in the complex plane by $D$, and given a disc $\Delta$, we denote the radius of $\Delta$ by $r(\Delta)$.

\bthm\label{main-theorem-plus}
For each $r>0$, there exists a sequence of open discs $\{D_k\}_{k=1}^\infty$ such that $\sum_{k=1}^\infty r(D_k)<r$ and such that setting $K=\ol D\setminus \bigcup_{k=1}^\infty D_k$, the uniform algebra $R(K)$ is nontrivial and strongly regular.
\ethm

A set of the form $K=\ol D\setminus \bigcup_{k=1}^\infty D_k$ with $\{D_k\}_{k=1}^\infty$ a sequence of open discs such that $\sum_{k=1}^\infty r(D_k)<\infty$ is called a \emph{Swiss cheese}.  (Sometimes in the literature a more restrictive definition of Swiss cheese is used.)  Thus the set $K$ in Theorem~\ref{main-theorem-plus} is a Swiss cheese.  For such a set $K$, the following standard result gives a very useful criterion insuring that $R(K)$ is nontrivial \cite[Lemma~24.1]{S1}.

\bthm\label{Stout}
Suppose that $\{D_k\}_{k=1}^\infty$ is a sequence of open discs in the complex plane such that $\sum_{k=1}^\infty r(D_k)<1$, and set $K=\ol D\setminus \bigcup_{k=1}^\infty D_k$.  Then $R(K)\not=C(K)$.
\ethm

The first example of a nontrivial {\it normal\/} uniform algebra was given by Robert McKissick \cite{McK}.  His example is $R(K)$ for a certain Swiss cheese $K$. Theorem~\ref{main-theorem-plus} is thus a sharpening of McKissick's result.

For $K$ a compact planar set such that $R(K)$ is nontrivial, the set of nonpeak points has positive planar measure~\cite[Theorem~26.8]{S1}.  Thus Theorem~\ref{main-theorem} gives an example of a strongly regular uniform algebra, on a metrizable space, with uncountably many nonpeak points.  Furthermore, replacing $R(K)$ by its restriction to its essential set (see the beginning of the proof of Theorem~\ref{bru-non-amenable}) the theorem gives an example with a dense set of nonpeak points.  All previously known examples of strongly regular uniform algebras on metrizable spaces had at most finitely many nonpeak points, and those on nonmetrizable spaces had at most finitely many points that were not generalized peak points \cite{F1}.  A strongly regular uniform algebra can never have nonzero \emph{bounded} point derivations.  The example in Theorem~\ref{main-theorem} is the first strongly regular uniform algebra known to have \emph{unbounded} point derivations.  (For $K$ a compact planar set, $R(K)$ has a nonzero point derivation at every nonpeak point \cite[Corollary~3.3.11]{Browder}).
Theorem~\ref{main-theorem} also gives the first example of a strongly regular uniform algebra with an infinite, in fact uncountable, Gleason part, since for $R(K)$ the Gleason part of every nonpeak point has positive planar measure \cite[Corollary~26.13]{S1}.  
An example with a two-point Gleason part was given by Feinstein \cite{F3}, and the same argument yields an example with an $n$-point Gleason part for any $n\in\Z_+$ \cite{GI}.
In addition, the example in Theorem~\ref{main-theorem} is the first strongly regular uniform algebra known to be finitely generated; for $K$ a compact planar set, $R(K)$ is always generated by two functions \cite[Corollary~24.4]{S1}.

In   \cite[Theorem~1.1]{FI},
Feinstein and the author introduced a general method for constructing essential uniform algebras.  Using this method they constructed an essential, natural, regular uniform algebra on the closed unit disc $\ol{D}$ \cite[Theorem~1.2]{FI}.  Repeating the proof of \cite[Theorem~1.2]{FI} with the strongly regular uniform algebra of Theorem~1.2 above in place of McKissick's normal uniform algebra shows that the result can be strengthened by replacing regularity by strong regularity.  (By a result of Wilken \cite[Lemma]{Wilken2}, every strongly regular uniform algebra is natural, so we omit mention of naturality from the statement of the theorem.)

\bthm
There exists an essential, strongly regular uniform algebra on the closed unit disc $\ol{D}$.
\ethm

A uniform algebra $A$ is \emph{weakly amenable} if there are no nonzero bounded derivations from $A$ into any commutative $A$-bimodule.  As proved in \cite[Theorem~1.5]{BCD}, weak amenability of $A$ is equivalent to the statement that there are no nonzero bounded derivations from $A$ into the dual commutative bimodule $A^*$.  It is standard that every trivial uniform algebra is weakly amenable (see \cite{Ka} or \cite[43.12]{BD}).

In \cite{FH}, Feinstein and Heath raised many questions including the following. 

\begin{question}
\cite[Question~5.5]{FH}
Is there a uniform algebra that is strongly regular but is not weakly amenable?
\end{question}

\begin{question}
\cite[Question~5.4]{FH}
Is there a uniform algebra that has bounded relative units but is not weakly amenable?
\end{question}

\begin{question}
\cite[Question~5.1]{FH}
Is there a nontrivial weakly amenable uniform algebra?
\end{question}

We will answer the first two of these questions in the affirmative.  The third question remains open.

\bthm\label{strongly-regular-non-amenable}
There exists a compact set $K$ in the complex plane such that $R(K)$ is strongly regular but not weakly amenable.
\ethm

\bthm\label{bru-non-amenable}
There exists an essential uniform algebra $A$ on a compact metrizable space such that $A$ has bounded relative units, but $A$ is not weakly amenable.
\ethm

To put these two theorems and their proofs in context, we recall some earlier examples.
Feinstein \cite{F3} constructed a compact planar set $K$ such that $R(K)$ has no nonzero bounded point derivations but $R(K)$ is not weakly amenable, and he constructed a uniform algebra $A$ on a compact metrizable space such that every point of the maximal ideal space of $A$ is a peak point for $A$ but $A$ is not weakly amenable.  Heath \cite{Heath} showed that the constructions could be modified so as to obtain examples that are regular. In \cite{FH}, Feinstein and Heath went further constructing an essential, regular uniform algebra on a compact metrizable space such that every point of the maximal ideal space $X$ of $A$ is a peak point while $A$ is not weakly amenable, and in addition, $A$ has bounded relative units at every point of a dense open subset of $X$.
Our approach to proving Theorems~\ref{strongly-regular-non-amenable} and~\ref{bru-non-amenable} is essentially to combine ideas from our proof of Theorem~\ref{main-theorem} with the methods used by Feinstein and Heath to construct 
their earlier examples.

We now describe our approach to proving Theorem~\ref{main-theorem-plus}.  Given a point $x\in\ol{D}$ and $r>0$, using McKissick's lemma \cite[Lemma~2]{McK} (see also 
\cite[Lemma~27.6]{S1} and \cite[Lemma~1.2]{Ko}), one can choose a sequence of open discs $\{D_k\}$ with $\sum r(D_k) < r$ and such that $R(\ol D\setminus \bigcup D_k)$ is strongly regular at $x$.  By repeating this at a countable set of points, Donald Chalice
obtained an $R(K)$ that is strongly regular at a countable dense set of points~\cite[pp.~302--303]{Chalice}.  However, to achieve strong regularity at every point more seems to be needed.  Specifically, if for a countable set $\{x_n\}$, we are to obtain the desired strong regularity from approximation of functions in $M_{x_n}$ by functions in $J_{x_n}$ for every $n\in \Z_+$, then we need some control over how large a disc about $x_n$ is contained in the zero set of an approximating function.  We will give a precise, general criterion (Lemma~\ref{criterion}) on a compact planar set $K$ that, given $x\in K$ and $s\in \Z_+$, insures that $\ol{J_x}\supset M_x^s$.  One might hope to apply McKissick's lemma to construct a compact set satisfying the criterion at every point $x\in K$ with $s=1$ and thereby obtain strong regularity.  However, the author was unable to do so.

A proof of McKissick's lemma, simpler than the original one, was given by Thomas K\"orner \cite{Ko}.  To get the control needed to satisfy the criterion in Lemma~\ref{criterion}, we will strengthen McKissick's lemma by giving certain uniform bounds on the functions used by K\"orner.  A surprising twist then arises in our argument.  Our bounds do not seem to enable us to directly satisfy the criterion of Lemma~\ref{criterion} with $s=1$.  Instead we obtain a Swiss cheese satisfying the criterion with $s=2$ (Theorem~\ref{half-way}).  Thus rather than a strongly regular uniform algebra, we obtain one in which $\ol{J_x}\supset M_x^2$ for every point $x$.  Finally, to conclude the proof of Theorem~\ref{main-theorem-plus}, we combine our construction with John Wermer's construction \cite[Theorem~1]{W1} of a Swiss cheese for which $\ol{M_x^2} =M_x$ for every $x$.

The condition that $\ol{J_x}\supset M_x^2$ for every $x\in K$ in the conclusion of Theorem~\ref{half-way} is actually equivalent to the condition that $\ol{J_x} = \ol{M_x^2}$ for every $x\in K$.
We digress from the main purpose of the paper to discuss this issue and in the process present some easy results of independent interest.

\bthm\label{ideal-space}
Fix $s\in \Z_+$.  If a uniform algebra $A$ on a compact space $K$ satisfies $\ol{J_x}\supset M_x^s$ for every $x\in K$, then $A$ is natural.
\ethm

\bcor\label{normal}
Fix $s\in \Z_+$.  If a uniform algebra $A$ on a compact space $K$ satisfies $\ol{J_x}\supset M_x^s$ for every $x\in K$, then $A$ is normal.
\ecor

We give here a simple proof of 
Theorem~\ref{ideal-space} by repeating an argument Raymond Mortini gave~\cite[Proposition~2.4]{Mortini} to reprove Wilken's result~\cite[Lemma]{Wilken2} that every strongly regular uniform algebra is natural.  Assume to get a contradiction that there exists a (nonzero) multiplicative linear function $\varphi$ on $A$ that is not evaluation at a point of $K$.  The hypothesis that $\ol{J_x}\supset M_x^s$ for every $x\in K$ implies that for each $x\in K$ there exists a function $f_x$ in $J_x$ such that $\varphi(f_x)=1$.  Then by the compactness of $K$, there is some finite subcollection $\{f_{x_1},\ldots, f_{x_n}\}$ of these functions whose zero sets cover $K$.  But then the product $f_{x_1} \cdots f_{x_n}$ is the zero function while $\varphi(f_{x_1} \cdots f_{x_n})=1$, a contradiction.  (Theorem~\ref{ideal-space} can also be proven by repeating Wilken's proof of \cite[Lemma]{Wilken2}.)
Corollary~\ref{normal} follows from Theorem~\ref{ideal-space} since the hypothesis on $A$ is easily seen to imply that $A$ is regular on $K$, and (as mentioned earlier) every uniform algebra that is regular on its maximal ideal space is normal \cite[Theorem~27.2]{S1}.

Let $A$ be a natural uniform algebra, and let $\mathfrak a$ be an ideal in $A$.  Let $E$ denote the hull of $\mathfrak a$ (the common zero set of the functions in $\mathfrak a$).  Set 
\[
J(E)=\{ f\in A:  \hbox{$f^{-1}(0)$ contains a neighborhood of $E$}\}.
\]
The ideal $\mathfrak a$ is said to be \emph{local} if $\mathfrak a\supset J(E)$.  The following result is standard \cite[Proposition~4.1.20(iv)]{D}.

\bthm\label{Garth}
Every ideal in a normal uniform algebra is local.
\ethm

As an immediate consequence of Corollary~\ref{normal} and Theorem~\ref{Garth} we get the following result, and in particular, we get that in Theorem~\ref{half-way} $\ol{J_x}=\ol{M_x^2}$ for every $x\in K$.

\bcor\label{equal}
Fix $s\in \Z_+$.  If a uniform algebra $A$ on a compact space $K$ satisfies $\ol{J_x}\supset M_x^s$ for every $x\in K$, then $\ol{J_x}= \ol{M_x^s}$ for every $x\in K$.
\ecor

Given that we will obtain our set $K$ such that $R(K)$ is strongly regular by separately imposing the conditions that $\ol{M_x^2}=M_x$ for every $x\in K$ and that $\ol {J_x}\supset M_x^2$ for every $x\in K$, it is natural to ask whether, for an $R(K)$ ($K$ compact planar), either of these conditions 
implies the other, and hence whether one of these conditions by itself is sufficient to insure that $R(K)$ is strongly regular.  The author recently gave examples showing that the first condition does not imply the second \cite[Theorems~1.4 and~1.5]{I}.  Conversely, it seems likely that the second condition does not imply the first.  In fact, we make the following conjecture.

\begin{conjecture}
For each integer $s\geq 2$, there exists a compact set $K\subset\C$ such that in $R(K)$ we have $\ol{J_x} \supset M_x^s$ for every $x\in K$ but there is some $y\in K$
such that $\ol{J_y}\not\supset M_y^{s-1}$.  Note that then by Corollary~\ref{equal} $\ol{J_x} = \ol{M_x^s}=\ol{M_x^{s+1}} = \ol{M_x^{s+2}} = \cdots$ for every $x\in K$ and $M_y\supsetneq \ol{M_y^2} \supsetneq \cdots \supsetneq \ol{M_y^s} = \ol{J_y}$.
\end{conjecture}

In the next section we define some terminology and notation already used above.  
In Section~\ref{bounds-on-Korner} we establish our uniform bounds on K\"orner's functions.  These bounds seem likely to have further applications.  In Section~\ref{main-construction} we give our general criterion, on a compact planar set $K$, for $R(K)$ to satisfy $\ol{J_x}\supset M_x^s$ for $x\in K$, and we construct a Swiss cheese $K$ for which the criterion holds for every $x\in K$ with $s=2$.
In Section~\ref{final-proof} we present two lemmas showing that certain inclusions of ideals in a uniform algebra $R(K)$ persist when we pass to a compact subset of $K$, and we conclude the proof of Theorem~\ref{main-theorem-plus} giving the existence of a nontrivial strongly regular $R(K)$.  The uniform algebras that are not weakly amenable in Theorems~\ref{strongly-regular-non-amenable} and~\ref{bru-non-amenable} are constructed in Sections~\ref{non-amenable} and~\ref{bounded-relative-units}, respectively.

%
%
%
%

\section{Terminology and Notation}\label{notation}

Those readers well versed in uniform algebra concepts may wish to skim this section and refer back to it as needed.

It is to be understood that all sequences, unions, and sums involving an index extend from $1$ to $\infty$; thus for instance $\{D_k\}$ means $\{D_k\}_{k=1}^\infty$, and $\bigcup D_k$ means $\bigcup_{k=1}^\infty D_k$.  
If $f$ is a function whose domain contains a subset $L$, we denote the restriction of $f$ to $L$ by $f|L$, and if $A$ is a collection of such functions, we denote the collection of restrictions of functions in $A$ to $L$ by $A|L$.
We denote the interior of a set $N$ by $N^\circ$.
The set of positive integers will be denoted by $\Z_+$.  Given a positive real number $c$, we denote by $c\ol D$ the closed disc of radius $c$ centered at the origin in the complex plane.

Recall that given a disc $\Delta$, we denote the radius of $\Delta$ by $r(\Delta)$.  
We will denote the distance from $\Delta$ to the coordinate axes $\R\cup i\R$ by $s_0(\Delta)$.  Explicitly, $s_0(\Delta)=\inf\{ |z-\zeta|: z\in \Delta,\ \zeta\in \R\cup i\R\}$.
More generally, given a point $a\in\C$, we will denote the distance from $\Delta$ to the union of the horizontal and vertical lines through $a$ by $s_a(\Delta)$.
We will denote the distance from $\Delta$ to the boundary $\partial I^2$ of the closed unit square $I^2=[0,1]\times[0,1]$ by $s(\Delta)$.

Throughout the paper all spaces will tacitly be required to be Hausdorff.  Let $X$ be a compact space.  We denote by $C(X)$ the algebra of all continuous complex-valued functions on $X$ with the supremum norm $\|f\|_\infty=\|f\|_X=\sup\{|f(x)|: x\in X\}$.  A \emph{uniform algebra} on $X$ is a closed subalgebra of $C(X)$ that contains the constants and separates the points of $X$.  A uniform algebra $A$ on $X$ is said to be 
\itemskip
\begin{enumerate}
\item[(a)] \emph{nontrivial} if $A\neq C(X)$,
\itemskip
\item[(b)] \emph{essential} if there is no proper closed subset of $E$ of $X$ such that $A$ contains every continuous complex-valued function on $X$ that vanishes on $E$,
\itemskip
\item[(c)] \emph{natural} if the maximal ideal space of $A$ is $X$ (under the usual identification of a point of $X$ with the corresponding multiplicative linear functional),
\itemskip
\item[(d)] \emph{regular on $X$} if for each closed set $K_0$ of $X$ and each point $x$ of $X\setminus K_0$, there exists a function $f$ in $A$ such that $f(x)=1$ and $f=0$ on $K_0$,
\itemskip
\item[(e)] \emph{normal on $X$} if for each pair of disjoint closed sets $K_0$ and $K_1$ of $X$, there exists a function $f$ in $A$ such that $f=1$ on $K_1$ and $f=0$ on $K_0$,
\end{enumerate}
The uniform algebra $A$ on $X$ is \emph{regular} or \emph{normal} if $A$ is natural and is regular on $X$ or normal on $X$, respectively.
In fact, every regular uniform algebra is normal \cite[Theorem~27.2]{S1}.
Also, if a uniform algebra $A$ is normal on $X$, then $A$ is necessarily natural \cite[Theorem~27.3]{S1}.  A uniform algebra that is regular on a compact space $X$ but is not natural was found by Hoffman and Singer~\cite{H-S} (or see \cite[pp.~187--190]{Hoffman}).

For $\varphi$ a multiplicative linear functional on the uniform algebra $A$, a \emph{point derivation} on $A$ at $\varphi$ is a linear functional $d$ on $A$ such that 
\begin{equation*}
d(fg)=d(f)\varphi(g) + \varphi(f)d(g)\qquad \hbox{for all $f$ and $g$ in $A$.}
\end{equation*}

A \emph{derivation} from $A$ to an $A$-bimodule $M$ is a linear map $D: A\ra M$ such that 
\begin{equation*}
D(fg)=D(f) \cdot g + f \cdot D(g)\qquad \hbox{for all $f$ and $g$ in $A$.}
\end{equation*}
A (point) derivation is said to be \emph{bounded} if it is bounded (i.e., continuous) as a linear map.

The notions of \emph{strongly regular} uniform algebra, \emph{weakly amenable} uniform algebra, and the ideals $M_x$ and $J_x$ were defined in the introduction.
When it is necessary to indicate with respect to which algebra the ideals are taken,
we will denote the ideals $M_x$ and $J_x$ in the uniform algebra $A$ by $M_x(A)$ and $J_x(A)$, respectively.

The \emph{essential set} $E$ for a uniform
algebra $A$ on $X$ is the unique smallest closed subset $E$ of $X$ such
that $A$ contains every continuous function on $X$ that vanishes on
$E$.   For a proof of the existence of the essential set and other details, see 
\cite[pp.~144--147]{Browder}.  Note that $A$ is essential if and only if the essential set for $A$ is $X$.

The point $x$ is said to be a \emph{peak point} for $A$ if there is a function $f$ in $A$ such that $f(x)=1$ and $|f(y)|<1$ for every $y\in X\sm \{x\}$.  The point $x$ is said to be a \emph{generalized peak point} if for every neighborhood $U$ of $x$ there exists a function $f$ in $A$ such that $f(x)=\|f\|=1$ and $|f(y)|<1$ for every $y\in X\sm U$.  When the space $X$ is metrizable, the notions of peak point and generalized peak point coincide.

It is well known that
for $K$ a compact planar set, the uniform algebra $R(K)$, which was defined in the introduction, is always natural~\cite[Theorem~24.5]{S1}.  It will be convenient to use the notation $R(K)$ not only when $K$ is a compact planar set, but more generally whenever $K$ is a compact subset of the Riemann sphere $\Rsphere$.  Thus for $K$ a compact subset of the Riemann sphere, $R(K)$ will denote the uniform closure on $K$ of the holomorphic rational functions with poles off $K$.

%
%
%
%

\section{Bounds on the functions of K\"orner}\label{bounds-on-Korner}

As mentioned in the introduction, we will use the functions that were introduced by K\"orner \cite{Ko} and used by him to simplify the proof of a lemma in McKissick's construction \cite{McK} of a normal uniform algebra.  In addition to the properties of these functions established by K\"orner, we will need certain uniform bounds on the functions.  Actually the functions we will use are not quite the same ones used by K\"orner.  Here is our modification of \cite[Lemma~1.2]{Ko}.

\blem\label{bounds}
There exist a sequence of rational functions $\{f_n\}_{n=1}^\infty$ and a constant $C_1>0$ such that for every $0<\vep < 1$ there is a sequence of open discs $\{\Dk\}_{k=1}^\infty$ in the plane such that 
\itemskip
\begin{enumerate}
\item[(a)] $\sum_{k=1}^\infty r(\Dk) \leq \vep$.
\itemskip
\item[(b)] The poles of the $f_n$ lie in $\unionDk$.
\itemskip
\item[(c)] The sequence $\{f_n\}$ converges uniformly on $\C\setminus \bigcup_{k=1}^\infty \Dk$ to a function $f=f_\vep$ that is identically zero outside $D$ and zero free in $D\setminus \unionDk$.
\itemskip
\item[(d)] $\unionDk\subset \{z: 1/2<|z|<1\}$.
\itemskip
\item[(e)] The discs $\Delta_1$, $\Delta_2$, $\ldots$ are disjoint.
\itemskip
\item[(f)] $f(0)=1$.
\itemskip
\item[(g)] $\|f_\vep\|_\infty\leq C_1\,\vep^{-1}$.
\end{enumerate}
Furthermore, given $M\in\Z_+$, the sequence $\{f_n\}$ can be chosen such that the derivatives $f^{(j)}(0)$ vanish for all $j=1,\ldots, M$.
\elem

In \cite{Ko}, the sequence $\{f_n\}$ depends on $\vep$, and this enables K\"orner to arrange to have in place of condition (d) that $\bigcup \Delta_k\subset \{z:1-\vep<|z|< 1\}$.  In the lemma above, the sequence $\{f_n\}$ is independent of $\vep$.  This has the advantage that any two limit functions $f_\vep$ and $f_{\vep'}$ agree where both are defined.  In particular, $f=f_\vep$ is independent of $\vep$ on $\{z: |z|<1/2\}$.  Another advantage is that this approach seems to yield a better dependence on $\vep$ of the bound in condition (g).

To obtain the uniform algebras that fail to be weakly amenable in Theorems~\ref{strongly-regular-non-amenable} and~\ref{bru-non-amenable}  we will need to show that condition (a) above can be replaced by a more stringent condition.  This will be done at the end of the section (see Lemma~ \ref{stronger-bounds}) so that the argument can be skipped by readers not concerned with failure of weak amenability.

The proof of Lemma~\ref{bounds} proceeds via several lemmas.  The following is \cite[Lemma~2.1]{Ko}.

\blem\label{Korner1}
If $N\geq 2$ is an integer and $h_N(z)=1/(1-z^N)$ then
\itemskip
\begin{enumerate}
\item[(i)] $|h_N(z)|\leq 2|z|^{-N}$\quad for $|z|^N\geq 2$,
\itemskip
\item[(ii)] $|1-h_N(z)|\leq 2|z|^N$\quad for $|z|^N\leq 2^{-1}$,
\itemskip
\item[(iii)] $h_N(z)\not= 0$\quad for all $z$.
\end{enumerate}
Furthermore, if $(8\log N)^{-1}>\delta>0$ then
\begin{enumerate}
\item[(iv)] $|h_N(z)|\leq 2\delta^{-1}$ provided only that $|z-w|\geq \delta N^{-1}$ whenever $w^N=1$.
\end{enumerate}
\elem

The following is \cite[Lemma~2.2]{Ko} except for the introduction of the scaling factor $\alpha$ in part (iv) and the addition of part (v).

\blem\label{Korner2}
If in Lemma~\ref{Korner1} we set $N=n2^{2n}$ with $n$ sufficiently large then for all $0<\alpha\leq 1$ we have
\begin{enumerate}
\item[(i)] $|h_N(z)|\leq (n+1)^{-4}$\quad for $|z|\geq 1+2^{-(2n+1)}$,
\itemskip
\item[(ii)] $|1-h_N(z)|\leq (n+1)^{-4}$\quad for $|z|\leq 1 - 2^{-(2n+1)}$,
\itemskip
\item[(iii)] $h_N(z)\not= 0$\quad for all $z$,
\itemskip
\item[(iv)] $|h_N(z)|\leq 2n^2\alpha^{-1}$ provided only that $|z-w|\geq n^{-3}2^{-2n}\alpha$ whenever $w^N=1$,
\itemskip
\item[(v)] $|h_N(z)|\leq 2^{2n+1}n^4\alpha^{-1}$ provided only that $|z-w|\geq n^{-5}2^{-4n}\alpha$ whenever $w^N=1$.
\end{enumerate}
\elem

\bpf
The proof is essentially a repetition of the proof of \cite[Lemma~2.2]{Ko}.
Since $\lim_{m\ra\infty} (1+ m^{-1})^m = e$, there is an integer $m_0$ such that $(1+\nobreak m^{-1})^m \geq e/2$ for all $m\geq m_0$.  (In fact, this inequality holds for all $m\in\Z_+$.) Thus for $n$ sufficiently large that $2n+1\geq m_0$ and $(e/2)^{n/2} \geq 2(n+1)^4$,
the condition $|z| \geq 1+2^{-(2n+1)}$ implies that 
$|z|^N \geq (1+2^{-(2n+1)})^N \geq (e/2)^{n/2} \geq 2(n+1)^4$.  Thus (i) follows from Lemma~\ref{Korner1}(i).  A similar argument yields (ii).
Parts (iv) and (v) follow from Lemma~\ref{Korner1}(iv) on setting $\delta= n^{-2}\alpha$ and $\delta= 2^{-2n}n^{-4}\alpha$, respectively, provided $n$ is chosen sufficiently large that $8\log N = 8(\log n + 2n \log 2) < n^2$.
\epf

The following is \cite[Lemma~2.3]{Ko} except for the introduction of the scaling factor $\alpha$.  

\blem\label{Korner3}
There exists an integer $n_0$ such that for each $n\geq n_0$ and each $0<\alpha \leq 1$ 
there exists a finite collection $\aa(n)$ of disjoint open discs and a rational function $g_n$ such that 
\begin{enumerate}
\item[(i)] $\sum\limits_{\Delta\in \aa(n)} r(\Delta)=n^{-2}\alpha$,
\itemskip
\item[(ii)] the poles of $g_n$ lie in $\bigcup\limits_{\Delta\in \aa(n)} \Delta$,
\itemskip
\item[(iii)] $|g_n(z)|\leq (n+1)^{-4}$\quad for $|z|\geq 1-2^{-(2n+1)}$,
\itemskip
\item[(iv)] $|1-g_n(z)|\leq (n+1)^{-4}$\quad for $|z|\leq 1 - 2^{-(2n-1)}$,
\itemskip
\item[(v)] $|g_n(z)|\leq 2n^2\alpha^{-1}$\quad for $z\notin \bigcup\limits_{\Delta\in\aa(n)} \Delta$,
\itemskip
\item[(vi)] $g_n(z)\not= 0$\quad for all $z$,
\itemskip
\item[(vii)] $\bigcup\limits_{\Delta\in\aa(n)} \Delta \subset \{z: 1-2^{-(2n-1)} \leq |z| \leq 1-2^{-(2n+1)}\}$.
\end{enumerate}
\elem

\bpf
The proof is essentially a repetition of the proof of \cite[Lemma~2.3]{Ko}.  Let $N=n2^{2n}$, $\omega=\exp(2\pi i/N)$, and $g_n(z)=h_N\bigl((1-2^{-2n})^{-1}z\bigr)$.  (A typo occurs in the definition of $N$ in \cite{Ko}.) Let $A_\alpha(n)$ be the collection of discs with radii $n^{-3}2^{-2n}\alpha$ and centers $(1-2^{-2n})\omega^j$ for $j=0,1,\ldots, N-1$.  Then the required results are either trivial or follow from Lemma~\ref{Korner2} after scaling $z$ by a factor of $1-2^{-2n}$ and making further routine estimates.
\epf

\bpf[Proof of Lemma~\ref{bounds}]
Let $n_0$ be as Lemma~\ref{Korner3}, and choose $m$ sufficiently large that $m\geq n_0$ and $2^{2m}>M$.   For each $n\geq m$ let 
$f_n=\prod_{j=m}^n g_j$, let $\alpha=\vep (\sum_{j=m}^\infty j^{-2})^{-1}$, and let $\{\Delta_k\}$ be an enumeration of the discs of $\bigcup_{j=m}^\infty \aa(j)$.
We may restrict attention to $0<\vep\leq \sum_{j=m}^\infty j^{-2}$ so that $0<\alpha\leq 1$.  Then conditions (a), (b), (d), and (e) are easily verified.

Set $C=\prod\limits_{j=1}^\infty \bigl[ 1+(j+1)^{-4} \bigr]$.
A tedious computation (or in K\"orner's words \lq\lq a simple induction\rq\rq) shows that for $z\notin \bigcup_{k=1}^\infty \Delta_k$ we have
\begin{align}
|f_n(z)| &\leq 2m^2C\alpha^{-1}   &{\rm if}\ |z|\leq 1 - 2^{-(2n-1)}\hphantom{\hskip 76.3pt} \nonumber\\
|f_n(z)|& \leq n^{-2}  \alpha^{-1} \leq 2m^2C\alpha^{-1}\quad   &{\rm if}\ 1-2^{-(2n-1)} \leq |z| \leq 1 - 2^{-(2n+1)} \nonumber\\
|f_n(z)|&\leq (n+1)^{-4} &{\rm if}\ 1- 2^{-(2n+1)} \leq |z|.\hphantom{\hskip 73.3pt} \nonumber
\end{align}
Condition (c) can now be proven as in \cite{Ko}.  Condition (f) is now easily verified.

For condition (g), note that we have from above that $|f_n(z)|\leq 2m^2C\alpha^{-1}$ for all $z\notin \bigcup_{k=1}^\infty \Delta_k$, and so the same inequality holds with $f$ in place of $f_n$.  Thus condition (g) holds with $C_1=2m^2C(\sum_{j=m}^\infty j^{-2})$.

The final assertion of the Lemma is evident since each $f_n$ is a function of $z^{2^{2m}}$ (and $m$ was chosen to satisfy $2^{2m}>M$).
\epf

The rest of this section can be skipped by readers not concerned with weak amenability.

Recall that given a disc $\Delta$, we denote the distance from $\Delta$ to the coordinate axes $\R\cup i\R$ by $s_0(\Delta)$.

\blem\label{stronger-bounds}
In Lemma~\ref{bounds} we can replace condition (a) that $\sum_{k=1}^\infty r(\Delta_k) \leq\vep$ by
\begin{enumerate}
\item[(a$'$)] $\sum_{k=1}^\infty r(\Delta_k)/s_0(\Delta_k)^2 \leq\vep$.
\end{enumerate}
\elem

Note that condition (a$'$) is stronger than condition (a) since $s_0(\Delta_k)\leq 1$ for all $k$.

The following is \cite[Lemma~2.4]{Heath} except for the introduction of the scaling factor $\alpha$.

\blem\label{stronger-Korner3}
In Lemma~\ref{Korner3} we can replace conditions (i) and (v) by 
\begin{enumerate}
\item[(i$'$)] $\sum\limits_{\Delta\in \aa(n)} r(\Delta)/s_0(\Delta)^2\leq n^{-2}\alpha$,
\itemskip
\item[(v$'$)] $|g_n(z)|\leq 2^{2n+1}n^4\alpha^{-1}$\quad for $z\notin \bigcup\limits_{\Delta\in\aa(n)} \Delta$.
\end{enumerate}
\elem

\bpf
The proof is essentially a repetition of the proof of \cite[Lemma~2.4]{Heath}; we repeat part of the argument for clarity.  Let $N=n2^{2n}$, $\omega=\exp(2\pi i/N)$, $\omega^{1/2}=\exp(\pi i/N)$, and $g=h_N\bigl((\omega^{1/2})^{-1} (1-2^{-2n})^{-1}z\bigr)$.  Let $A_\alpha(n)$ to be the collection of discs with radii $n^{-5}2^{-4n}\alpha$ and centers $(1-2^{-2n})\omega^{1/2}\omega^j$ for $j=0,1,\ldots, N-1$.  Then, with the exception of (i$'$), the results are either trivial or follow from Lemma~\ref{Korner2} on \lq\lq scaling" by a factor of $\omega^{1/2}(1-2^{-2n})$.  For the proof of (i$'$) we refer the reader to \cite{Heath} noting that of course one must now invoke Lemma~\ref{Korner2} above in place of \cite[Lemma~2.3]{Heath}.  (The argument in \cite{Heath} contains the erroneous inequality 
\[
\frac{(r+1/2)\pi}{n2^{2n+2}} - n^{-5} 2^{-4n} \geq \left(\frac{1}{2}\right) \frac{2r+1}{n2^{2n}},
\]
but this does not matter since it is true that (for all $n$ sufficiently large) the left hand side is greater than some fixed constant times the right hand side, and that suffices for the argument.)
\epf

\bpf[Proof of Lemma~\ref{stronger-bounds}]
The proof is similar to the proof of Lemma~\ref{bounds}.  Choose $m$ as in the proof of Lemma~\ref{bounds}.  Let $f_n=\prod_{j=m}^n g_j$, let $\alpha=\vep (\sum_{j=m}^\infty j^{-2})^{-1}$, and let $\{\Delta_k\}$ be an enumeration of the discs of $\bigcup_{j=m}^\infty \aa(j)$ with $\aa(j)$ as in Lemma~\ref{stronger-Korner3}.
We may restrict attention to $0<\vep\leq \sum_{j=m}^\infty j^{-2}$ so that $0<\alpha\leq 1$.  Then conditions (a), (b), (d), and (e) are easily verified.

Set $C=\prod\limits_{j=1}^\infty \bigl[ 1+(j+1)^{-4} \bigr]$.
A tedious computation shows that for $z\notin \bigcup_{k=1}^\infty \Delta_k$ we have
\begin{align}
|f_n(z)| &\leq 2^{2m+1}m^4C\alpha^{-1}   &{\rm if}\ |z|\leq 1 - 2^{-(2n-1)}\hphantom{\hskip 76.5pt} \nonumber\\
|f_n(z)|& \leq 2^{-2n} 2^{4m+5}  \alpha^{-1}  &{\rm if}\ 1-2^{-(2n-1)} \leq |z| \leq 1 - 2^{-(2n+1)} \nonumber\\
|f_n(z)|&\leq [(n+1)!]^{-4}(m!)^4\quad &{\rm if}\ 1- 2^{-(2n+1)} \leq |z|.\hphantom{\hskip 73.3pt}\nonumber
\end{align}
For notational convenience set $\tilde C=\max\{2^{2m+1}m^4C,\ 2^{-2n} 2^{4m+5}\}$, so that these inequalities become
\begin{align}
|f_n(z)| &\leq \tilde C\alpha^{-1}   &{\rm if}\ |z|\leq 1 - 2^{-(2n-1)}\hphantom{\hskip 76.5pt}
\nonumber\\
|f_n(z)|& \leq \tilde C \alpha^{-1}  &{\rm if}\ 1-2^{-(2n-1)} \leq |z| \leq 1 - 2^{-(2n+1)} \nonumber\\
|f_n(z)|&\leq [(n+1)!]^{-4}(m!)^4\quad &{\rm if}\ 1- 2^{-(2n+1)} \leq |z|.\hphantom{\hskip 73.5pt}\nonumber
\end{align}

As in K\"orner \cite{Ko}, apply the trivial equality
$$|f_{n+1}(z) - f_n(z)| = |f_n(z)|\, |1-g_{n+1}(z)|.$$
This gives
\begin{align}
|f_{n+1}(z) - f_n(z)|&\leq \tilde C\alpha^{-1} (n+2)^{-4}   
\nonumber\\
&\phantom{\leq \tilde C\alpha^{-1} (n+1)^{-4}}{\rm if}\ |z|\leq 1 - 2^{-(2n+1)}\nonumber\\
\nonumber\\
|f_{n+1}(z) - f_n(z)|& \leq [(n+1)!]^{-4} (m!)^4 [1+2^{2n+3}(n+1)^4\alpha^{-1}]\nonumber\\
&\phantom{\,[(n+1)!]^{-4} (m!)^4}{\rm if}\ 1- 2^{-(2n+1)} \leq |z|. \nonumber
\end{align}
Consequently, for all $n$ sufficiently large, we have $|f_{n+1}(z) - f_n(z)|\leq \tilde C\alpha^{-1} n^{-4}$ for all $z\notin \bigcup \Delta_k$.  Therefore, by the Weierstrass $M$-test, $(f_n)$ converges uniformly to some function $f$ on $\C\setminus \bigcup \Dk$.

That $f(z)=0$ for $|z| \geq 1$ is evident.
The proof that $f(z)\not = 0$ for $|z|<1$ is the same as in K\"orner \cite{Ko}.  
Thus condition (c) is verified.
The remainder of the proof is essentially the same as in the proof of Lemma~\ref{bounds}.
\epf

%
%

\section{The main construction}\label{main-construction}

Our goal in this section is to prove the existence of a Swiss cheese $K$ such that in the uniform algebra $R(K)$ we have $\ol {J_x}\supset M_x^2$ for every $x\in K$.  We will begin with a lemma that provides a general criterion, given $s\in \Z_+$, for the inclusion $\ol {J_x}\supset M_x^s$ to hold.  We will then prove a technical lemma that will enable us to satisfy that criterion at every point with $s=2$.  Finally using the two lemmas, we will construct the desired Swiss cheese.

\blem\label{criterion}
Let $K\subset\C$ be compact, let $s\in \Z_+$, and let $x\in K$.  Suppose that for every $\sigma>0$ and $\eta>0$ there is an open disc $\Delta$ that contains $x$ and has radius $r(\Delta)\leq \sigma$, and such that, denoting the center of $\Delta$ by $a$, there is a function $g$ in $R(K)$ such that $g$ is identically zero on $\Delta\cap K$ and $\| (z-a)^s - g \|_K < \eta$.  Then in the uniform algebra $R(K)$, we have $\ol {J_x} \supset {M^s_x}$.
\elem

\bpf
Since the closed ideal generated by the function $(z-x)^s$ is $\ol{M^s_x}$, it suffices to show that $(z-x)^s$ is in $\ol {J_x}$.  Let $\sigma>0$ and $\eta>0$ be arbitrary.  Let $\Delta$, $a$, and $g$ be as in the statement of the lemma.  
Let $m=\sup_{z\in K} \bigl| (d/dz) (z^s) \bigr|$.  Then
\[
\| (z-x)^s - (z-a)^s \|_K \leq m \sigma.
\]
Thus
\begin{align}
\| (z-x)^s - g \|_K&\leq \| (z-x)^s - (z-a)^s \|_K + \| (z-a)^s - g \|_K \nonumber\\
&<m\sigma + \eta.\nonumber
\end{align}
Since $g$ is in $J_x$, and $\sigma>0$ and $\eta>0$ are arbitrary, this shows that $(z-x)^s$ is in $\ol {J_x}$, as desired.
\epf

Recall that the restriction to $(1/2)\ol D$ of the function $f=f_\vep$ in Lemma~\ref{bounds} is independent of $\vep$.  Thus, in particular, $\|f''\|_{(1/4)\ol D}$ is a well-defined number.

\blem\label{technical}
Let $C_1$ and $f=f_\vep$ be as in Lemma~\ref{bounds}.  Set 
\[
C=\min\{ 2^{1/2}\|f''\|_{(1/4)\ol D}^{-1/2}\,\,,\,\, 2^{-5/2} C_1^{-1/2} \}.
\]
Then for 
$0<\vep\leq 1$ and $\eta>0$ and 
\[
\sigma=\sigma(\vep,\eta)=C\eta^{1/2} \vep^{1/2}
\]
there is a sequence of disjoint open discs $\{D_k\}_{k=1}^\infty$ contained in the annulus $\{z: \sigma < |z| < 2\sigma\}$ such that
\begin{enumerate}
\item[(a)] $\sum_{k=1}^\infty r(D_k) < 4\sigma\vep$, and
\item[(b)] there is a function $h\in R(\Rsphere\setminus \bigcup_{k=1}^\infty D_k)$ such that $h=0$ on $\{|z|\leq \sigma\}$ and such that for all $z\in \C \setminus \bigcup_{k=1}^\infty D_k$ we have 
\[
|z^2 - z^2 h(z) |\leq\eta.
\]
\end{enumerate}
\elem

\bpf
Given $0<\vep\leq 1$, choose a sequence of disjoint open discs $\{\Delta_k\}$ as in Lemma~\ref{bounds}.  Set $M=\|f''\|_{(1/4)\ol D}$.  By the final assertion of Lemma~\ref{bounds} we may assume that $f'(0)=0$.  Then,
\begin{equation}\label{F}
| f_\vep (w) - 1| = | f_\vep (w) - f_\vep(0) | \leq (M/2) |w|^2 \quad \hbox{for}\ |w| \leq 1/4.
\end{equation}

Set $h(z)=h_\vep(z)=f_\vep(\sigma/z)$ (and $h(0)=0$).  Let $D_k$ be the image of $\Delta_k$ under the map $z\mapsto \sigma/z$.  Since each $\Delta_k$ lies in $\{z: 1/2 <|z|<1\}$, we have $\sum r(D_k) \leq 4\sigma \sum r(\Delta_k) < 4\sigma\vep$.  The function $h$ is defined on $\C\setminus \bigcup D_k$ and is a uniform limit there of rational functions with poles in $\bigcup D_k$.  Furthermore, $h=0$ on $\{ |z|\leq \sigma\}$.  Also inequality (\ref{F}) gives
\begin{equation}\label{first}
|h(z) -1 | \leq (M/2) \sigma^2 |z^{-2}| \quad \hbox{for}\ |z|\geq 4\sigma.
\end{equation}

It remains to be shown that 
\begin{equation}
|z^2 - z^2 h(z) | = |z^2|\, |h(z) -1 | \leq \eta \quad \hbox{for all}\ z\in \C \setminus {\textstyle \bigcup D_k}.\nonumber
\end{equation}
For $|z|\geq 4\sigma$, we have by inequality (\ref{first})
\[
|z^2|\, |h(z)-1|\leq (M/2)\sigma^2 = (M/2)C^2\eta\vep\leq \eta\vep\leq \eta.
\]
For $|z|\leq 4\sigma$ (and $z\notin \bigcup D_k$), applying condition (g) of Lemma~\ref{bounds} yields
\begin{align}
|z^2|\, |h(z)-1| &\leq (4\sigma)^2 \bigl(\|h\|_\infty +1\bigr) \nonumber\\
& =  (4\sigma)^2 \bigl(\|f_\vep\|_\infty +1\bigr) \nonumber\\
& \leq (4\sigma)^2 \bigl(C_1\vep^{-1} +1\bigr) \nonumber\\
& \leq (4\sigma)^2 \bigl(2C_1\vep^{-1}\bigr) \nonumber\\
&\leq\eta\nonumber.
\end{align}
The lemma is proved.
\epf

A simple translation argument yields the following immediate corollary.

\bcor\label{translation}
Let $\sigma=\sigma(\vep,\eta)$ be as in Lemma~\ref{technical}.  
Given $a\in \C$, there is a sequence of disjoint open discs $\{D_k\}_{k=1}^\infty$ contained in the annulus $\{z: \sigma < |z-a| < 2\sigma\}$ such that
\begin{enumerate}
\item[(a)] $\sum_{k=1}^\infty r(D_k) < 4\sigma\vep$, and
\itemskip
\item[(b)] there is a function $h\in R(\Rsphere\setminus \bigcup_{k=1}^\infty D_k)$ such that $h=0$ on $\{|z-a|\leq \sigma\}$ and such that for all $z\in \C \setminus \bigcup_{k=1}^\infty D_k$ we have 
\[
\bigl|(z-a)^2 - (z-a)^2 h(z)\bigr|\leq\eta.
\]
\end{enumerate}
\ecor

\begin{observation}\label{observation}
Given $c>0$, there is a number $M\leq 9/c^2$ such that the square $[-1,1]\times[-1,1]$, and hence the disc $\ol D$, can be covered by $M$ open discs of radius $c$.  To see this, note that $[-1,1]\times[-1,1]\supset \ol D$ can be expressed as the union of $\bigl( \lceil 2/c \rceil \bigr)^2 \leq 9/c^2$ squares of side length $c$, and each such square is contained in the open disc of radius $c$ with center the center of the square.
\end{observation}

We can now prove the main result of this section.

\bthm\label{half-way}
For each $r>0$, there exists a sequence of open discs $\{D_k\}_{k=1}^\infty$ such that $\sum_{k=1}^\infty r(D_k)<r$ and such that setting $K=\ol D\setminus \bigcup_{k=1}^\infty D_k$ we have, in the uniform algebra $R(K)$, that $\ol {J_x} \supset M_x^2$ for every $x\in K$.
\ethm

In fact, as discussed in the introduction, the algebra $R(K)$ in the theorem actually satisfies $\ol {J_x} = \ol{M_x^2}$ for every $x\in K$.

\bpf
Fix $r>0$.  Let $\sigma=\sigma(\vep, \eta)$ and $C$ and $C_1$ be as in Lemma~\ref{technical}.  Note that
\[
\vep/\sigma = C^{-1}\eta^{-1/2} \vep^{1/2}. 
\]
It follows trivially that setting $\eta_n=1/n$, choosing a sequence $\{\vep_n\}$ going to zero fast enough, and setting $\sigma_n=\sigma(\vep_n,\eta_n)$, we can arrange to have $36\sum \vep_n/\sigma_n < r$ and $\sigma_n\ra0$ as $n\ra\infty$.

By Observation~\ref{observation}, for each $n=1,2,\ldots$, we can cover $\ol D$ by a collection $\sd_n$ of $M_n\leq 9/\sigma_n^2$ open discs of radius $\sigma_n$.
Given a disc $\Delta$ in $\sd_n$, let $a$ denote the center of $\Delta$ so that $\Delta=\{|z-a|<\sigma_n \}$.  Apply
Corollary~\ref{translation} to choose a sequence of open discs $\{D_k^a\}$ contained in the annulus $\{z: \sigma_n < |z-a| < 2\sigma_n\}$ such that $\sum r(D_k^a) < 4\sigma_n\vep_n$ and there is a function $h\in R(\Rsphere\setminus \bigcup D_k^a)$ such that $h=0$ on $\{|z-a|\leq\sigma_n\}$ and such that for all $z\in\C\sm\bigcup D_k^a$ we have $\bigl|(z-a)^2 - (z-a)^2 h(z)\bigr|\leq\eta_n$.  Do this for each disc in each $\sd_n$.  Let $\{D_k\}$ be an enumeration of all the discs so chosen.  Then $\sum r(D_k) < 36 \sum \vep_n/\sigma_n <r$. 

Set $K=\ol D\setminus\bigcup D_k$.  Let $x\in K$ be arbitrary.  Given $\sigma>0$ and $\eta>0$, choose $n\in\Z_+$ large enough that $\sigma_n<\sigma$ and $\eta_n<\eta$.  There is a disc $\Delta$ in $\sd_n$ such that $x$ is in $\Delta$.  Then $r(\Delta)=\sigma_n<\sigma$, and letting $a$ denote the center of $\Delta$, there exists a function $h\in R(\Rsphere\sm \bigcup D_k)$ that is identically zero on $\Delta\cap K$ and satisfies $\| (z-a)^2 - (z-a)^2 h \|_K \leq \eta_n<\eta$.  Therefore, applying Lemma~\ref{criterion} shows that $\ol {J_x}\supset M_x^2$, as desired.
\epf

%
%
%
%

\section{Proof of the main theorem}\label{final-proof}

In this section we complete the proof of Theorem~\ref{main-theorem-plus} by combining Theorem~\ref{half-way} with Wermer's theorem \cite{W1} on the existence of a Swiss cheese $K$ such that the uniform algebra $R(K)$ has no nonzero bounded point derivations.
It is well known (and easy to show) \cite[p.~64]{Browder} that for a uniform algebra $A$ on a compact space $K$ and a point $x\in K$, there exists a bounded point derivation on $A$ at $x$ if and only if $M_x^2$ is not dense in $M_x$.  Thus Wermer's result can be restated as follows.

\bthm\label{Wermer}
For each $r>0$, there exists a sequence of open discs $\{D_k\}_{k=1}^\infty$ such that $\sum_{k=1}^\infty r(D_k)<r$ and such that setting $K=\ol D\setminus \bigcup_{k=1}^\infty D_k$ we have, in the uniform algebra $R(K)$, that $\ol{M_x^2}=M_x$ for every $x\in K$.
\ethm

We will need two simple lemmas showing that certain inclusions of ideals in a uniform algebra $R(K)$ persist when we pass to a compact subset $L$ of $K$. 

\blem\label{subsetJ}
Given compact sets $L\subset K\subset \C$, given $s\in \Z_+$, and given a point $x\in L$, if 
$\ol{J_x(R(K))} \supset M_x(R(K))^s$, then 
$\ol{J_x(R(L))} \supset M_x(R(L))^s$.
\elem

\bpf
Since the closed ideal generated by the function $(z-x)^s$ is $\ol{M_x(R(L))^s}$, it suffices to show that $(z-x)^s$ is in $\ol{J_x(R(L))}$.    
Fix $\vep > 0$ arbitrary.  By hypothesis, $(z-x)^s$ is in $\ol {J_x(R(K))}$.  Thus there is a function $f$ in $R(K)$ that vanishes on a neighborhood of $x$ in $K$ such that $\|f-(z-x)^s\|_K <\vep$.  Then the function $f|L$ is in $R(L)$, vanishes on a neighborhood of $x$ in $L$, and satisfies $\| f-(z-x)^s\|_L< \vep$.  Consequently, $(z-x)^s$ is in $\ol{J_x(R(L))}$, as desired.
\epf

\blem\label{subset2}
Given compact sets $L\subset K\subset \C$ and given a point $x\in L$, if 
$\ol{M_x(R(K))^2}=M_x(R(K))$, then 
$\ol{M_x(R(L))^2}=M_x(R(L))$.
\elem

\bpf
This can be proven by an argument analogous to the one just given for Lemma~\ref{subsetJ}.  Alternatively, the lemma follows immediately from Hallstrom's criterion \cite[Theorem~1]{Hall} for the existence of a nonzero bounded point derivation on $R(K)$. 
\epf

We can now finish the proof of Theorem~\ref{main-theorem-plus}.  

\bpf[Proof of Theorem~\ref{main-theorem-plus}]
Without loss of generality assume that $r<1$.  By Theorem~\ref{half-way}, there exists a sequence of open discs $\{\Delta_k^I\}_{k=1}^\infty$ such that $\sum r(\Delta_k^I)< r/2$ and such that setting $K_1=\ol D\setminus \bigcup \Delta_k^I$ we have $\ol{J_x(R(K_1))} \supset M_x(R(K_1))^2$ for every $x\in K_1$.  By Theorem~\ref{Wermer}, 
there exists a sequence of open discs $\{\Delta_k^W\}_{k=1}^\infty$ such that $\sum r(\Delta_k^W)< r/2$ and such that setting $K_2=\ol D\setminus \bigcup \Delta_k^W$ we have that $\ol{M_x(R(K_2))^2}=M_x(R(K_2))$ for every $x\in K_2$.  Let $\{D_k\}$ be an enumeration of the collection of discs $\{\Delta_k^I\} \cup \{\Delta_k^W\}$.  Set $K=K_1\cap K_2=\ol D \setminus \bigcup D_k$.  Then the uniform algebra $R(K)$ is nontrivial by Theorem~\ref{Stout}.  Furthermore, 
Lemmas~\ref{subsetJ} and~\ref{subset2} yield that $\ol {J_x(R(K))} \supset \ol{M_x(R(K))^2} = M_x(R(K))$ for every $x\in K$.  Thus $R(K)$ is strongly regular.
\epf

%
%
%
%

\section{Strong regularity without weak amenability}\label{non-amenable}

In this section we establish Theorem~\ref{strongly-regular-non-amenable}.  The construction of the desired compact set $K$ is similar to the construction of the set $K$ in Theorem~\ref{main-theorem-plus} but somewhat more elaborate.
Our approach to proving that our uniform algebra is not weakly amenable is the same as in the papers \cite{F3, FH, Heath} of Feinstein and Heath.  Given a compact planar set $K$ denote by $R_0(K)$ the set of restrictions to $K$ of the holomorphic rational functions with poles off $K$.  Suppose that $\mu$ is a measure on $K$ such that the bilinear form defined on $R_0(K)\times R_0(K)$ by
$$(f,g)\mapsto \int_K f'g\, d\mu$$
is bounded.  Then as noted in \cite{F0}, we can extend the form by continuity to $R(K)\times R(K)$ and obtain a bounded derivation $D:R(K)\ra R(K)^*$ such that for $f$ and $g$ in $R_0(K)$,
$$(Df)(g)= \int_K f'g\, d\mu.$$
Such a derivation is the zero derivation if and only if $\mu$ annihilates $R(K)$.

We will prove the following theorem.

\bthm\label{thm-bounded-derivation}
For each $C>0$ there exists a compact planar set $K$ obtained by deleting from the closed unit square $I^2$ a countable union of open discs such that the boundary $\partial I^2$ of $I^2$ is contained in the essential set for $R(K)$, $R(K)$ is strongly regular, but for all $f$ and $g$ in $R_0(K)$,
\begin{equation}\label{bounded-derivation}
\left | \int_{\partial I^2} f'(z) g(z)\, dz \right | \leq C \|f\|_K \|g\|_K.
\end{equation}
\ethm

Theorem~\ref{strongly-regular-non-amenable} is an immediate consequence, for taking $K$ to be as in Theorem~\ref{thm-bounded-derivation} and $a$ to be a point of $I^2\sm K$, note that $\int_{\partial I^2} 1/(z-a)\, dz = 2\pi i \not= 0$, so the above discussion shows that there is a nonzero bounded derivation $D:R(K)\ra R(K)^*$ such that for $f$ and $g$ in $R_0(K)$,
$$(Df)(g)=\int_{\partial I^2} f'(z)g(z)\, dz.$$

\begin{remark}\label{remark}
Given a uniform algebra $A$ on a compact space $X$ and a closed subset $E$ of $X$, there is an obvious algebraic isomorphism of the restriction algebra $A|E$ with the quotient algebra $A/I(E)$, where $I(E)$ denotes the ideal of functions in $A$ vanishing on $E$.  Thus $A|E$ can be regarded as a Banach algebra using the quotient norm on $A/I(E)$.
Applying an argument of Feinstein \cite[p.~2393]{F3} then shows that, with $K$ the compact planar set of the above theorem, any uniform algebra with a restriction isomorphic to $R(K)|\partial I^2$ must fail to be weakly amenable.
\end{remark}

The proof of Theorem~\ref{thm-bounded-derivation} will use several preliminary lemmas.  The first of these is a
minor modification of \cite[Lemma~2.1]{F3} of Feinstein.  Its proof is essentially identical to that of Feinstein's lemma and hence omitted.

Recall the notations $r(\Delta)$, $s_0(\Delta)$, $s_a(\Delta)$, and $s(\Delta)$ introduced in the introduction.

\blem\label{Feinstein}
Let $\{D_k\}_{k=1}^\infty$ be a sequence of open discs in the complex plane whose closures are contained in the interior of the unit square $I^2$.  Set $K=I^2\sm \bigcup_{k=1}^\infty D_k$.  Let $f$ and $g$ be in $R_0(K)$.  Then
$$\left | \int_{\partial I^2} f'(z) g(z) \, dz \right | \leq 4\pi \|f\|_K \|g\|_K \sum_{k=1}^\infty \frac{r(D_k)}{s(D_k)^2}.$$
\elem

\blem\label{essential-set}
Let $\{D_k\}_{k=1}^\infty$ be a sequence of open discs in the complex plane whose closures are contained in the interior of the unit square $I^2$.  Set $K=I^2\sm \bigcup_{k=1}^\infty D_k$. Suppose that $\sum_{k=1}^\infty r(D_k)/ s(D_k)<\infty$.  Then $\partial I^2$ is contained in the essential set for $R(K)$.
\elem

\bpf
For convenience set $r_k=r(D_k)$ and $s_k= s(D_k)$.
Choose $N$ such that $\sum_{k=N}^\infty r_k/s_k < 1/2$.  Set $\delta=\min\{s_1,\ldots, s_{N-1}\}$.  Let $E$ be an arbitrary closed disc of radius $r(E)<\delta/2$ that is contained in $I^2$ and intersects $\partial I^2$. We will show that $R(E\sm \bigcup D_k)\not= C(E\sm \bigcup D_k)$.  Since every (relatively) open set of $I^2$ that intersects $\partial I$ contains such a disc $E$, the lemma follows.

By Theorem~\ref{Stout} and a trivial scaling argument, to show that $R(E\sm \bigcup D_k)\not= C(E\sm \bigcup D_k)$ it suffices to show that $\sum\limits_{\{k: D_k\cap E\not= \emptyset\} }r_k < r(E)$.  Now note that
\begin{align}
\sum\limits_{\{k: D_k\cap E\not= \emptyset\} }r_k/r(E)&\leq \sum_{\{k: s_k\leq 2r(E)\} } r_k/r(E)\nonumber\\
&\leq \sum_{\{k: s_k\leq 2r(E)\} } 2r_k/s_k \nonumber\\
& \leq 2 \sum_{k=N}^\infty r_k/ s_k \nonumber\\
&<1\nonumber.
\end{align}
\vskip -16pt
\epf

\blem\label{Hallstrom}
Let $K$ be a compact set in the complex plane, and let $x\in K$.  If there is a neighborhood $U$ of $x$ in $K$ such that $R(\ol U)$ has no nonzero bounded point derivations at $x$, then $R(K)$ also has no nonzero bounded point derivations at $x$
\elem

\bpf
This is immediate from Hallstrom's criterion \cite[Theorem~1]{Hall} for the existence of a nonzero bounded point derivation on $R(K)$. 
\epf

\bthm\label{Wermer-not-amenable}
For each $r>0$, there exists a sequence of open discs $\{D_k\}_{k=1}^\infty$ such that $\sum_{k=1}^\infty r(D_k)/ s(D_k)^2<r$ and such that setting $K=I^2\setminus \bigcup_{k=1}^\infty D_k$ we have, in the uniform algebra $R(K)$, that $\ol{M_x^2}=M_x$ for every $x\in K$.
\ethm

\bpf
Choose a countable collection $\{U_j\}$ of open discs that covers $(I^2)^\circ$ with  the closure of each $U_j$ is contained in $(I^2)^\circ$.   By Theorem~\ref{Wermer}, we can choose, for each $j$, a sequence of open discs $\{\Delta_l^j\}_{l=1}^\infty$ with each $\Delta_l^j$ contained in $U_j$ such that $\sum_{l=1}^\infty r(\Delta_l^j) < r s(U_j)^2 /2^j$ and such that in the uniform algebra $R(\ol U_j \sm \bigcup_{l=1}^\infty\Delta_l^j)$ we have $\ol{M_x^2} = M_x$ for every $x\in \ol U_j \sm \bigcup_{l=1}^\infty \Delta_l^j$.  
Let $\{D_k\}$ be an enumeration of the collection $\{ \Delta_l^j: j=1, 2, \ldots, l=1, 2, \ldots\}$.  Then $\sum r(D_k)/ s(D_k)^2<r$.  Set $K=I^2\setminus \bigcup D_k$.  Then the equality 
$\ol{M_x^2}=M_x$ holds for $x\in K \cap (I^2)^\circ$ by 
Lemma~\ref{Hallstrom} and holds for $x\in \partial I^2$ because each point of $\partial I^2$ is a peak point for $R(K)$.
\epf

Recall that the restriction to $(1/2)\ol D$ of the function $f=f_\vep$ in Lemmas~\ref{bounds} and~\ref{stronger-bounds} is independent of $\vep$.  Thus, in particular, $\|f'''\|_{(1/4)\ol D}$ is a well-defined number.

\blem\label{technical-non-amenable}
Let $C_1$ and $f=f_\vep$ be as in Lemma~\ref{stronger-bounds} {\rm (}see also Lemma~\ref{bounds}{\rm )}.  Set 
\[
C=\min\{6^{1/3}\|f'''\|_{(1/4)\ol D}^{-1/3}\,\,,\,\, 2^{-7/3} C_1^{-1/3} \}.
\]
Then for 
$0<\vep\leq 1$ and $\eta>0$ and 
\[
\rho=\rho(\vep,\eta)=C\eta^{1/3} \vep^{1/3}
\]
there is a sequence of disjoint open discs $\{D_k\}_{k=1}^\infty$ contained in the annulus $\{z: \rho < |z| < 2\rho\}$ such that 
\begin{enumerate}
\item[(a)] $\sum_{k=1}^\infty r(D_k)/s_0(D_k)^2 < 4\vep/\rho$, and
\item[(b)] there is a function $h\in R(\Rsphere\setminus \bigcup_{k=1}^\infty D_k)$ such that $h=0$ on $\{|z|\leq \rho\}$ and such that for all $z\in \C \setminus \bigcup_{k=1}^\infty D_k$ we have 
\[
|z^3 - z^3 h(z)|\leq\eta.
\]
\end{enumerate}
\elem

\bpf
The proof is similar to the proof of Lemma~\ref{technical}.
Given $0<\vep\leq 1$, choose a sequence of disjoint open discs $\{\Delta_k\}$ as in Lemma~\ref{bounds} such that condition (a$'$) of Lemma~\ref{stronger-bounds} holds.  Set $M=\|f'''\|_{(1/4)\ol D}$.  By the final assertion of Lemma~\ref{bounds} we may assume that $f''(0)=0$.  Then,
\begin{equation}\label{FF}
| f_\vep (w) - 1| = | f_\vep (w) - f_\vep(0) | \leq (M/6) |w|^3 \quad \hbox{for}\ |w| \leq 1/4.
\end{equation}

Set $h(z)=h_\vep(z)=f_\vep(\rho/z)$ (and $h(0)=0$).  Let $D_k$ be the image of $\Delta_k$ under the map $z\mapsto \rho/z$.  Since each $\Delta_k$ lies in $\{z: 1/2 <|z|<1\}$, we have $r(D_k) \leq 4\rho r(\Delta_k)$.  From the geometry of the conformal map $z\mapsto \rho/z$ we see that $s(D_k) \geq \rho s(\Delta_k)$.  Consequently, 
$\sum r(D_k)/s_0(D_k)^2 \leq 4 \sum r(\Delta_k)/\rho s_0(\Delta_k)^2 < 4\vep/\rho$.  The function $h$ is defined on $\C\setminus \bigcup D_k$ and is a uniform limit there of rational functions with poles in $\bigcup D_k$.  Furthermore, $h=0$ on $\{ |z|\leq \rho\}$.  Also inequality (\ref{FF}) gives
\begin{equation}\label{first2}
|h(z) -1 | \leq (M/6) \rho^3 |z^{-3}| \quad \hbox{for}\ |z|\geq 4\rho.
\end{equation}

It remains to be shown that 
\begin{equation}
|z^3 - z^3 h(z)| = |z^3|\, |h(z) -1 | \leq \eta \quad \hbox{for all}\ z\in \C \setminus {\textstyle \bigcup D_k}.\nonumber
\end{equation}
For $|z|\geq 4\rho$, we have by inequality (\ref{first2})
\[
|z^3|\, |h(z)-1|\leq (M/6)\rho^3 = (M/6)C^3\eta\vep\leq \eta\vep\leq \eta.
\]
For $|z|\leq 4\rho$ (and $z\notin \bigcup D_k$), applying condition (g) of Lemma~\ref{bounds} and Lemma~\ref{stronger-bounds} yields
\begin{align}
|z^3|\, |h(z)-1| &\leq 
(4\rho)^3 \bigl(\|f_\vep\|_\infty +1\bigr) \nonumber\\
& \leq (4\rho)^3 \bigl(2C_1\vep^{-1}\bigr) \nonumber\\
&\leq\eta\nonumber.
\end{align}
The lemma is proved.
\epf

A simple translation argument yields the following immediate corollary.

\bcor\label{translation-not-amenable}
Let $\rho=\rho(\vep,\eta)$ be as in Lemma~\ref{technical-non-amenable}.  
Let $a\in \C$ be fixed.  Then
there is a sequence of open discs $\{D_k\}_{k=1}^\infty$ contained in the annulus $\{z: \rho < |z-a| < 2\rho\}$ such that
\begin{enumerate}
\item[(a)] $\sum_{k=1}^\infty r(D_k)/ s_a(D_k)^2 < 4\vep/\rho$, and
\itemskip
\item[(b)] there is a function $h\in R(\Rsphere\setminus \bigcup_{k=1}^\infty D_k)$ such that $h=0$ on $\{|z-a|\leq \rho\}$ and such that for all $z\in \C \setminus \bigcup_{k=1}^\infty D_k$ we have 
\[
\bigl|(z-a)^3 - (z-a)^3 h(z)\bigr|\leq\eta.
\]
\end{enumerate}
\ecor

\bthm\label{half-way-not-amenable}
For each $r>0$, there exists a sequence of open discs $\{D_k\}_{k=1}^\infty$ contained in the open unit square $(I^2)^\circ$ such that $\sum_{k=1}^\infty r(D_k)/ s(D_k)^2<r$ and such that setting $K= I^2\setminus \bigcup_{k=1}^\infty D_k$ we have, in the uniform algebra $R(K)$, that $\ol {J_x} \supset M_x^3$ for every $x\in K$.
\ethm

In fact, as discussed in the introduction, the algebra $R(K)$ in the theorem actually satisfies $\ol {J_x} = \ol{M_x^3}$ for every $x\in K$.

\bpf
Fix $r>0$.  The square $I^2$ is the union of three (disjoint) sets: the interior of $I^2$, the boundary of $I^2$ minus the set of corners, and the set of corners.  We will work in turn on each of these three sets.

Let $\sigma=\sigma(\vep, \eta)$ and $C$ and $C_1$ be as in Lemma~\ref{technical}.  Note that
\[
\vep/\sigma = C^{-1}\eta^{-1/2} \vep^{1/2}. 
\]
It follows trivially that setting $\eta_n=1/n$, choosing a sequence $\{\vep_n\}$ going to zero fast enough, and setting $\sigma_n=\sigma(\vep_n,\eta_n)$, we can arrange to have $36\sum n^2\vep_n/\sigma_n < r/3$ and $\sigma_n\ra0$ as $n\ra\infty$.

For each $n\in \Z_+$ such that $\frac{1}{n} + 4\sigma_n<1$, let $Q_n$ be the square
$$Q_n=\Bigl[-1+\bigl(\textstyle\frac{1}{n} + 4\sigma_n\bigr), 1 - \bigl(\frac{1}{n} + 4\sigma_n\bigr)\Bigr]^2.$$
By Observation~\ref{observation}, we can cover $Q_n$ by a collection $\sd_n^1$ of  $M_n\leq 9/\sigma_n^2$ open discs of radius $\sigma_n$.  We may assume that each disc in $\sd_n^1$ intersects $Q_n$.

Given a disc $\Delta$ in $\sd_n^1$, let $a$ denote the center of $\Delta$ so that $\Delta=\{|z-a|<\sigma_n \}$.  Choose a sequence of open discs $\{D_k^a\}$ as in Corollary~\ref{translation} (with $\vep=\vep_n$, $\eta=\eta_n$, and $\sigma=\sigma_n$), and note that each disc $D_k^a$ chosen satisfies $s(D_k^a)>1/n$.  Carry out this procedure for each disc $\Delta$ in each $\sd_n^1$.  Let $\{D_k^1\}$ be an enumeration of all the discs so chosen.  
Then $\sum r(D_k^1)/s(D_k^1)^2 < \sum (9/\sigma_n^2)(4\sigma_n\vep_n) n^2 = 36 \sum n^2\vep_n/\sigma_n <r/3$. 

Now let $\rho=\rho(\vep, \eta)$ and $C$ and $C_1$ be as in Lemma~\ref{technical-non-amenable}.  Note that
\[
\vep/\rho^2 = C^{-2}\eta^{-2/3} \vep^{1/3}. 
\]
It follows trivially that setting $\eta_n=1/n$, choosing a sequence $\{\vep_n\}$ going to zero fast enough, and setting $\rho_n=\rho(\vep_n,\eta_n)$, we can arrange to have $8\sum \vep_n/\rho_n^2 < r/3$ and $\rho_n\ra0$ as $n\ra\infty$.

For each $n\in \Z_+$ such that $8\rho_n<1$, let $L_n$ be the union of the four line segments obtained by deleting from $\partial I^2$ the four discs of radius $8\rho_n$ whose centers are the four corners of $I^2$.
We can cover $L_n$ by a collection $\sd_n^2$ of  $N_n\leq 2/\rho_n$ open discs of radius $\rho_n$ with center in $L_n$.  

Given a disc $\Delta$ in $\sd_n^2$, let $a$ denote the center of $\Delta$ so that $\Delta=\{|z-a|<\rho_n \}$.  Choose a sequence of open discs $\{D_k^a\}$ as in Corollary~\ref{translation-not-amenable} (with $\vep=\vep_n$, $\eta=\eta_n$, and $\rho=\rho_n$), and note that each disc $D_k^a$ chosen satisfies $s(D_k^a)\geq s_a(D_k^a)$.  Carry out this procedure for each disc $\Delta$ in each $\sd_n^2$.  Let $\{D_k^2\}$ be an enumeration of all the discs so chosen.  
Then $\sum r(D_k^2)/s(D_k^2)^2 < \sum (2/\rho_n)(4\vep_n/\rho_n) = 8  \sum \vep_n/\rho_n^2 <r/3$. 

Again setting $\eta_n=1/n$, choose a new sequence $\{\vep_n\}$ such that again setting $\rho_n=\rho(\vep_n,\eta_n)$, we have $16 \sum \vep_n/\rho_n <r/3$ and $\rho_n\ra0$ as $n\ra\infty$.  For each $n\in\Z_+$, let $\sd_n^3$ be the collection whose members are the four discs of radius $\rho_n$ whose centers are the four corners of $I^2$.

Given a disc $\Delta$ in $\sd_n^3$, let $a$ denote the center of $\Delta$ so that $\Delta=\{|z-a|<\rho_n \}$.  Choose a sequence of open discs $\{D_k^a\}$ as in Corollary~\ref{translation-not-amenable} (with $\vep=\vep_n$, $\eta=\eta_n$, and $\rho=\rho_n$), and note that each disc $D_k^a$ chosen satisfies $s(D_k^a) = s_a(D_k^a)$.  Carry out this procedure for each disc $\Delta$ in each $\sd_n^3$.  Let $\{D_k^3\}$ be an enumeration of all the discs so chosen.  
Then $\sum r(D_k^3)/s(D_k^3)^2 < \sum 4(4\vep_n/\rho_n) = 16 \sum \vep_n/\rho_n <r/3$. 
 
Let $\{D_k\}$ be an enumeration of those discs in $\{D_k^1\} \cup \{D_k^2\} \cup \{D_k^3\}$ that are contained in $(I^2)^\circ$.  Then $\sum r(D_k)/s(D_k)^2 < r$.
Set $K=\ol D\setminus\bigcup D_k$.  
Essentially the same
reasoning used to conclude the proof of Theorem~\ref{half-way} now shows that it follows from Lemma~\ref{criterion} that 
$\ol {J_x}\supset M_x^3$ for every $x\in K$ (and in fact $\ol {J_x}\supset M_x^2$ for every $x\in K\cap (I^2)^\circ$).
\epf

\bpf[Proof of Theorem~\ref{thm-bounded-derivation}]
By Theorem~\ref{half-way-not-amenable}, there exists a sequence of open discs $\{\Delta_k^I\}$ such that $\sum r(\Delta_k^I)/ s(\Delta_k^I)^2< C/8\pi$ and such that setting $K_1= I^2\setminus \bigcup \Delta_k^I$ we have $\ol{J_x(R(K_1))} \supset M_x(R(K_1))^3$ for every $x\in K_1$.  By Theorem~\ref{Wermer-not-amenable}, 
there exists a sequence of open discs $\{\Delta_k^W\}$ such that $\sum r(\Delta_k^W)/ s(\Delta_k^W)^2< C/8\pi$ and such that setting $K_2= I^2\setminus \bigcup \Delta_k^W$ we have that $\ol{M_x(R(K_2))^2}=M_x(R(K_2))$ for every $x\in K_2$.  Let $\{D_k\}$ be an enumeration of the collection of discs $\{\Delta_k^I\} \cup \{\Delta_k^W\}$.  Set $K=K_1\cap K_2= I^2 \setminus \bigcup D_k$.  Then
Lemmas~\ref{subsetJ} and~\ref{subset2} yield that $R(K)$ is strongly regular.  (Note that the condition $\ol{M_x^2}=M_x$ implies $\ol{M_x^3}=M_x$.)  
Since $\sum r(D_k)/ s(D_k)^2 < C/4\pi$,
Lemma~\ref{Feinstein} shows that inequality~(\ref{bounded-derivation}) holds for all $f$ and $g$ in $R_0(K)$.  Lemma~\ref{essential-set} shows that $\partial I^2$ is contained in the essential set for $R(K)$.
\epf

%
%
%
%

\section{Bounded relative units without weak amenability}\label{bounded-relative-units}

In this section we prove Theorem~\ref{bru-non-amenable} by applying Brian Cole's method of root extensions to the uniform algebra given by Theorem~\ref{strongly-regular-non-amenable}.  We begin by recalling some aspects of Cole's construction \cite{Cole} (see also \cite[Section~19]{S1}).

Let $A$ be a uniform algebra on a compact space $X$, and let $\sf$ be a (nonempty) subset of $A$.  Endow $\C^\sf$ with the product topology.  Let $p_1:X\times \C^\sf\ra X$ and $p_f:\xcf\ra\C$ denote the projections given by $p_1(x, (z_g)_{g\in \sf})=x$ and $p_f(x, (z_g)_{g\in\sf})=z_f$.  Define $X_\sf\subset \xcf$ by
\[
X_\sf = \{\, y\in\xcf:\bigl(p_f(y)\bigr)^2=f\bigl(p_1(y)\bigr) \ \hbox{for all $f\in\sf$}\,\},
\]
and let $A_\sf$ be the uniform algebra on $X_\sf$ generated by the set of functions $\{\, f\circ p_1: f\in A\} \cup \{\, p_f: f\in\sf\}$.  On $X_\sf$ we have $p_f^2=f\circ p_1$ for every $f\in \sf$.  Set $\pi=p_1|X_\sf$, and note that $\pi$ is surjective.  There is an isometric embedding $\pi^*:A\ra A_\sf$ given by $\pi^*(f)=f\circ \pi$.  

We call the uniform algebra $A_\sf$ or the pair $(A_\sf, X_\sf)$, the $\sf$-extension of $A$, and we call $\pi$ the associated surjection.  Note that if $X$ is metrizable and $\sf$ is countable, then $X_\sf$ is metrizable also.   Given $x\in X$, if $\sf$ is contained in $M_x$, then the set 
$\pi^{-1}(x)$ consists of a single point.  

To prove 
Theorem~\ref{bru-non-amenable}, we will iterate 
the above extension process to obtain an infinite sequence of 
uniform algebras and then take a direct limit to obtain the desired uniform algebra.

We will need the following lemma of Feinstein and Heath~\cite[Lemma~4.3]{FH}.

\blem\label{bru}
Let $A$ be a uniform algebra on $X$ and $x\in X$.  Suppose that, for each compact subset $E$ of $X\sm \{x\}$, there exists a neighborhood $U$ of $x$ and a function $f\in A$ such that
\begin{enumerate}
\item[(i)] $f|U=1$,
\item[(ii)] $f|E=0$,
\item[(iii)] For each $k\in\Z_+$ there is a function $g\in A$ with $g^{2^k}=f$.
\end{enumerate}
Then $A$ has bounded relative units at $x$.
\elem

The next lemma, whose elementary proof we omit, is a modification of a lemma of Feinstein~\cite[Lemma~3.5]{F1}.

\blem\label{mod}
Let $A$ be a normal uniform algebra on a compact metrizable space $X$, and let $F$ be a closed subset of $X$.  Then there exists a countable subset $\sf$ of $A$ consisting of functions each vanishing identically on a neighborhood of $F$ such that for each point $x\in X\sm F$, and for each compact subset $E$ of $X\sm\{x\}$, there exists a neighborhood $U$ of $x$, and a function $f\in \sf$ such that $f|U=1$ and $f|E=0$.
\elem

\bpf[Proof of Theorem~\ref{bru-non-amenable}]
Let $K$ be the compact planar set given by Theorem~\ref{strongly-regular-non-amenable}.  Let $X_0$ be the essential set for $R(K)$.  Strong regularity of $R(K)|X_0$ follows trivially from strong regularity of $R(K)$.  Furthermore, $R(K)|X_0=R(X_0)$ (by \cite[Lemma~3.2.5]{Browder} for instance).  Thus $R(X_0)$ is essential and strongly regular.  Note also that $R(X_0)|\partial I^2=R(K)|\partial I^2$.

We will construct a sequence of uniform algebras $\{A_m\}_{m=0}^\infty$.  First set $A_0=R(X_0)$, and set $F_0=\partial I^2$.  By Lemma~\ref{mod}, there is a countable subset $\sf_0$ of $A_0$ consisting of functions each vanishing identically on a neighborhood of $F_0$ such that for each point $x\in  X_0\sm F_0$, and for each compact subset $E$ of $X_0\sm \{x\}$, there exists a neighborhood $U$ of $x$, and a function $f\in \sf_0$ such that $f|U=1$ and $f|E=0$.  Let $(A_1,X_1)$ be the $\sf_0$-extension of $A_0$, and let $\pi_1:X_1\rightarrow X_0$ be the associated surjection.  Because each member of $\sf_0$ is identically zero on $F_0$, the map $\pi_1$ takes $\pi_1^{-1}(F_0)$ homeomorphically onto $F_0$.  Let $F_1=\pi_1^{-1}(F_0)$.  By \cite[Theorem~2.4]{F1}, $A_1$ is normal.

We then iterate this process to obtain a sequence $\{A_m,X_m, \pi_m, F_m, \sf_m)\}_{m=0}^\infty$ where each $(A_m, X_m)$ is the $\sf_{m-1}$-extension of $(A_{m-1}, X_{m-1})$, each $\pi_m: X_m\rightarrow X_{m-1}$ is the associated surjection, $F_m=\pi_m^{-1}(F_{m-1})$, and $\sf_m$ is a countable subset of $A_m$ consisting of functions each vanishing identically on a neighborhood of $F_m$ such that for each function $f$ in $\sf_{m-1}$ the function $f\circ \pi_m$ is the square of a function in $\sf_m$ and such that for each point $x\in X_m\sm F_m$, and for each compact subset $E$ of $X_m\sm \{x\}$, there exists a neighborhood $U$ of $x$, and a function $f\in \sf_m$ such that $f|U=1$ and $f|E=0$.  Because each member of $\sf_{m-1}$ is identically zero on $F_{m-1}$, the map $\pi_m$ takes $F_m$ homeomorphically onto $F_{m-1}$.  Finally we take the direct limit of the system of uniform algebras $\{A_m\}$.  Explicitly, we set
$$X_\omega=\Bigl\{(x_j)_{j=0}^\infty \in \textstyle\prod\limits_{j=0}^\infty X_j : \pi_{m+1}(x_{m+1})=x_m
\hbox{\ for all } m=0, 1, 2, \ldots\Bigr\},$$
and letting $q_m:X_\omega\rightarrow X_m$ be the restriction of the canonical projection $\prod_{j=0}^\infty X_j\rightarrow X_m$, we let $A_\omega$ be the closure of $\bigcup_{m=0}^\infty \{h\circ q_m: h\in A_m\}$ in $C(X_\omega)$.  Set $F_\omega=\{(x_j)_{j=0}^\infty \in X_\omega : x_0\in F_0\}$.  Set $\pi=q_0$.  Then $\pi$ maps $F_\omega=\pi^{-1}(F_0)$ homeomorphically onto $F_0$.
By \cite[Corollary~2.9]{F1}, $A_\omega$ is normal.

Note that $X_\omega$ is metrizable.

To prove that $A_\omega$ is essential, we first prove by induction that $A_m$ is essential for each $m=0, 1, 2, \ldots$.  Recall that $A_0=R(X_0)$ is essential.  Assume as the induction hypothesis that $A_{m-1}$ is essential.  Denote by $\Z_2$ the multiplicative group consisting of $1$ and $-1$ and set $\sf=\sf_{m-1}$.  The product group $\Z_2^{\sf}$ acts as a topological transformation group on $X_m$ as follows: Given $\gamma\in \Z_2^\sf$ and $(x,z)\in X_m\subset X_{m-1}\times \C^\sf$, let $\gamma(x,z)$ be the point $(x,z')$ where the $f$th coordinate of $z'$ is $\gamma_f z_f$.  The uniform algebra $A_m$ is invariant under this action in the sense that the function $(x,z)\mapsto f\bigl(\gamma(x,z)\bigr)$ is in $A_m$ for each $\gamma\in \Z_2^\sf$ whenever $f$ is in $A_m$.  Consequently, the essential set for $A_m$ must also be invariant under the action of $Z_2^\sf$ on $X_m$.  Since $\Z_2^\sf$ acts transitively on each fiber of $\pi_m$, it follows that if $A_m$ is not essential, then there must be an open set $U$ of $X_{m-1}$ such that $\pi_m^{-1}(U)$ lies outside the essential set for $A_m$.  It is standard (see \cite[Theorem~19.1]{S1}) that given a uniform algebra $A$ on a compact space $X$, a subset $\sf$ of $A$, the $\sf$-extension $A_\sf$ of $A$ with associated surjection $\pi$, and a function $g\in C(X)$ such that the function $g\circ\pi$ is in $A_\sf$, then the function $g$ must be in $A$.  Consequently, the condition that $\pi_m^{-1}(U)$ lies outside the essential set for $A_m$ implies that $U$ lies outside the essential set for $A_{m-1}$, a contradiction.  Thus $A_m$ must be essential.

If $A_\omega$ is not essential, then there must be some $n\in \Z_+$ and some open set $W$ in $X_n$ such that 
$q_n^{-1}(W)$ lies outside
the essential set for $A_\omega$.  It is standard that if $f\in C(X_n)$ is such that $f\circ q_n$ is in $A_\omega$, then $f$ is in $A_n$.  (This can be derived from \cite[Lemma~19.3]{S1} noting that discarding the terms $A_0,\ldots, A_{n-1}$ from the system of uniform algebras $\{A_m\}$ does not affect the direct limit $A_\omega$.)  Consequently, the condition that $q_n^{-1}(W)$ lies outside the essential set for $A_\omega$ implies that $W$ lies outside the essential set for $A_n$, a contradiction.  Thus $A_\omega$ must be essential.

To prove that $A_\omega$ has bounded relative units, we consider separately points in the fibers of $\pi$ lying over $X_0\sm F_0$ and over $F_0$.  At points over $X_0\sm F_0$ we apply Lemma~\ref{bru}.  Consider a point $x\in X_\omega$ lying over $X_0\sm F_0$ and a compact subset $E$ of $X_\omega\sm\{x\}$.  For some $n\in \Z_+$, the point $q_n(x)$ is outside $q_n(E)$.  Then there is a neighborhood $V$ of $q_n(x)$ in $X_n$ and a function $h\in \sf_n$ such that $h|U=1$ and $h|(q_n(E))=0$.  Then the neighborhood $U=q_n^{-1}(V)$ and the function $f=h\circ q_n$ satisfy the conditions in Lemma~\ref{bru}.  Consequently, $A_\omega$ has bounded relative units at $x$.

Now consider a point $x$ lying over $F_0$.  To show that $A_\omega$ has bounded relative units at $x$ we first note that by a theorem of Feinstein \cite[Theorem~1.5]{F2} it is enough to show that $x$ is a peak point for $A_\omega$ and that $A_\omega$ is strongly regular at $x$.  That $x$ is a peak point is immediate since given a function $h$ in $A_0$ that peaks at $\pi(x)$, the function $h\circ \pi$ peaks at $x$, since $x$ is the only point in the fiber $\pi^{-1}(\pi(x))$.  Because the set $\bigcup_{m=0}^\infty \{h\circ q_m: h\in A_m\}$ is dense in $A_\omega$, to show that $A_\omega$ is strongly regular at $x$, it suffices to show that $A_m$ is strongly regular at $q_m(x)$ for each $m=0, 1, 2,\ldots$.  Recall that $A_0$ is strongly regular.  Assume for the purpose of induction that $A_{m-1}$ is strongly regular at $q_{m-1}(x)$.  Every function in $A_m$ that vanishes at $q_m(x)$ can be uniformly approximated by a function $f$ of the form 
$$f = \pi_m^*(f_0)+\sum_{u=1}^s \pi_m^*(f_u) g_u$$
where $f_0, f_1,\ldots,f_s\in A_{m-1}$, the function $f_0$ vanishes at   $q_{m-1}(x)$, and each $g_u$ is a product of functions of the form $p_f$ for $f\in \sf_{m-1}$.  Because every function in $\sf_{m-1}$ vanishes identically on a neighborhood of  $q_{m-1}(x)$, the sum $\sum_{u=1}^s \pi_m^*(f_u) g_u$ vanishes identically on a neighborhood of $q_m(x)$.  Because $A_{m-1}$ is strongly regular at $q_{m-1}(x)$, the function $f_0$ can be uniformly approximated by a function in $A_{m-1}$ vanishing identically on a neighborhood of $q_{m-1}(x)$, and hence $\pi_m^*(f_0)$ can be uniformly approximated by a function in $A_m$ vanishing identically in a neighborhood of $q_m(x)$.  Therefore, $A_m$ is strongly regular at $q_m(x)$.

Finally we show that $A_\omega$ is not weakly amenable.  There is a standard norm-decreasing linear operator $T:A_\omega\rightarrow A_0$ such that $T(f\circ \pi)=f$ for every $f\in A_\omega$.  (See \cite[Lemma~19.3]{S1} for instance.)  It is easily seen from the way this operator $T$ is defined that for each point $x$ of $X_0$ whose fiber $\pi^{-1}(x)$ consists of a single point, we have for each function $g\in A_\omega$ that $\bigl(T(g)\bigr)(x)= g(\pi^{-1}(x)\bigr)$.  (The operator $T$ is essentially given by averaging over the fibers of $\pi$.)  Consequently, for each $g\in A_\omega$ we have $g|F_\omega= \bigl(T(g)|F_0\bigl)\circ\pi$.  It follows that the restriction algebras $A_\omega|F_\omega$ and $A_0|F_0$ are isomorphic, and thus $A_\omega$ fails to be weakly amenable by Remark~\ref{remark}.
\epf

\section*{Acknowledgment}

The author thanks the referee for doing an exceptionally thorough job and making numerous observations and suggestions that improved the paper.


\begin{thebibliography}{BC87}

\bibitem{BCD}
W. G. Bade, P. C. Curtis, Jr., and H. G. Dales, Amenability and weak amenability for Beurling and Lipschitz algebras, {\it Proc.\ London Math.\ Soc.\/} 55 (1987), 359--377.

\bibitem{BD}
F. F. Bonsall and J. Duncan, {\it Complete normed algebras}, Springer, New York, 1973.

\bibitem{Browder} 
A. Browder, {\it Introduction to Function Algebras}, W. A. Benjamin, New York, NY, 1969.

\bibitem{Chalice}
D. R. Chalice, S-algebras on sets in $C^n$, {\it Proc.\ Amer.\ Math.\ Soc.\/} 39 (1973), 300--304.

\bibitem{Cole}
B.~J.~Cole, {\it One-Point Parts and the Peak Point Conjecture\/}, Ph.D. Thesis, Yale University, 1968.

\bibitem{D}
H.~G.~Dales, {\it Banach Algebras and Automatic Continuity}, London Mathematical Society Monographs, New Series, Volume 24,
The Clarendon Press, Oxford, 2000.

\bibitem{F0} J.\ F.\ Feinstein, Weak ($F$)-amenability of $R(X)$, {\it Proceedings of the Centre for Mathematical Analysis\/}, Australian National University, Canberra, 21 (1989), 97--125.

\bibitem{F1}
J.\ F.\ Feinstein,   A non-trivial, strongly regular uniform algebra, {\it J. London Math.\ Soc.\/} (2), 45 (1992), 288--300.

\bibitem{F2}  {J.\ F.\ Feinstein}, {Regularity conditions for a Banach function algebra}, in {\it Function Spaces\/} (Edwardsville, Il., 1994), Lecture Notes in Pure Appl.\ Math.\ 172, Dekker, New York, 1995, 117--122.

\bibitem{F3}
J.\ F.\ Feinstein,   A counterexample to a conjecture of S. E. Morris, {\it Proc.\ Amer.\ Math.\ Soc.\/} 132 (2004), 2389--2397.

\bibitem{FH} {J.\ F.\ Feinstein and M.\ J.\ Heath}, {Regularity and amenability conditions for uniform algebras\/}, Function spaces, 
{\it Contemporary Math.},  435 (2007),   159--169.

\bibitem{FI} J. F. Feinstein and A. J. Izzo, A general method for constructing essential uniform algebras, {\it Studia Math.} 246 (2019), 47--61.

\bibitem{GI} S. N. Ghosh and A. J. Izzo, One-point Gleason parts and point derivations in uniform algebras, {\it Studia Math.\/} 270 (2023), 323--337.

\bibitem{Hall}  A. P. Hallstrom, On bounded point derivations and analytic capacity, {\it J.\ Functional Analysis}, 4 (1969), 153--165.

\bibitem{Heath} M. J. Heath, A note on a construction of J. F. Feinstein, {\it Studia Math.\/} 169 (2005), 63--70.

\bibitem{Hoffman} 
K. Hoffman, {\it Banach Spaces of Analytic Functions}, Dover Publications Inc., Mineola, NY, 1962.

\bibitem{H-S} K. Hoffman and I. M. Singer, {Maximal algebras of continuous functions}, {\it Acta Math.\/} 103 (1960), 217--241.

\bibitem{I}
A. J. Izzo, {A normal uniform algebra that fails to be strongly regular at a peak point},
{\it Pacific J.\ Math.\/} 331 (2024), 77--97.

\bibitem{Ka}
H. Kamowitz, Cohomology groups of commutative Banach algebras, {\it Trans.\ Amer.\ Math.\ Soc.\/} 102 (1962), 352--372.

\bibitem{Ko} T. W. K\"orner, A cheaper Swiss cheese, {\it Studia Math.\/} 83 (1986), 33--36.

\bibitem{McK} 
R.~J.~McKissick,   A nontrivial normal sup norm algebra, {\it Bull.\  Amer.\ Math.\ Soc.\/} 69 (1963), 391--395. 

\bibitem{Mortini}
R. Mortini, Closed and prime ideals in the algebra of bounded analytic functions, {\it Bull.\ Austral.\ Math.\ Soc.\/} 35 (1987), 213--229.

\bibitem{S1}
 E.~L.~Stout, {\it The Theory of Uniform Algebras\/}, Bogden and 
Quigley, New York, 1971.

\bibitem{W1}
J.~Wermer,   Bounded point derivations on certain Banach algebras, {\it  J.\ Functional  Analysis},  1 (1967), 28--36.

\bibitem{Wilken1} 
D. R. Wilken, Approximate normality and function algebras on the interval and the circle, in: {\it Function Algebras\/}, F. Birtel (ed.),
Scott, Foresman and Co., Chicago, 1966, pp.~98--111.


\bibitem{Wilken2} 
D. R. Wilken, A note on strongly regular function algebras, {\it Canad.\ J. Math.}\ 21 (1969) 912--914.

\end{thebibliography}
\end{document}